\documentclass[12pt]{article}
\usepackage{eucal, amssymb,amsmath,graphicx, amsfonts, latexsym}

\begin{document} \parskip=5pt plus1pt minus1pt \parindent=0pt

\newcommand{\E}{{\rm E}}
\newcommand{\Var}{{\rm Var}}
\newcommand{\re}{{\rm e}}
\newcommand{\rg}{{\rm g}}
\newcommand{\rh}{{\rm h}}
\newcommand{\ru}{{\rm u}}
\newcommand{\rv}{{\rm v}}
\newcommand{\bn}{{\mbox{\boldmath $n$}}}
\newcommand{\ba}{{\mbox{\boldmath $a$}}}
\newcommand{\bA}{{\mbox{\boldmath $A$}}}
\newcommand{\btheta}{\bm{\theta}}

\title{Stochastic epidemic models: a survey}
\author{ Tom Britton, Stockholm University\thanks{Department of
Mathematics, Stockholm University, SE-106 91 Stockholm, Sweden. {\it
E-mail}: tom.britton@math.su.se}}
\date{\today}
\maketitle

\begin{abstract}
\noindent This paper is a survey paper on stochastic epidemic
models. A simple stochastic epidemic model is defined and exact and
asymptotic model properties (relying on a large community) are
presented. The purpose of modelling is illustrated by studying effects
of vaccination and also in terms of inference procedures for important
parameters, such as the basic reproduction number and the
critical vaccination coverage. Several generalizations towards
realism, e.g.\ multitype and household epidemic models,
are also presented, as is a model for endemic diseases.

\end{abstract}

{\bf Keywords}: Basic reproduction number, critical vaccination
coverage, household epidemic, multitype epidemic, stochastic epidemic, threshold theorem.


\section{Introduction}\label{intro}

Early modelling contributions for infectious disease spread were often
for specific diseases. For example Bernoulli (1760) aimed at
evaluating the effectiveness a certain technique of variolation
against smallpox, and Ross (1911) modelled the transmission of
malaria. One of the first more general and rigorous study was made by
Kermack and McKendrick (1927). Later important contributions were for
example by Bartlett (1949) and Kendall (1956), both also considering
stochastic models.

Early models were often deterministic and the type of questions that were
adressed were for example: Is it possible that there is a big
outbreak infecting a positive fraction of the community?, How many
will get infected if the epidemic takes off?, What are the effects of
vaccinating a given community fraction prior to the arrival of the
disease?, What is the endemic level? As problems were resolved,
 the
simple models were generalised in several ways towards making them
more realistic. Some such extensions were for example to allow for a
community where there are different \emph{types} of individual,
allowing for non-uniform mixing between individuals (i.e.\ infectious
individuals don't infect all individuals equally likely), for example
due to social or spatial aspects, and to allow seasonal variations.

Another generalisation of the initial simple deterministic epidemic
model was to study \emph{stochastic} epidemic models. 
A stochastic model is
of course preferable when studying a small community. But, even when
considering a large community, which deterministic models
primarily are aimed for,  some
additional questions can be raised when considering
stochastic epidemic models. For example: What is the
\emph{probability} of a major outbreak?, and for models describing an
endemic situation: How long is the disease likely to
persist (with or without intervention)? Later stochastic models have
also shown to be advantageous when the contact structure in the
community contains small complete graphs; households and other
local social networks being common examples. Needless to say, both
deterministic and stochastic epidemic models have their important
roles to play however, the focus in the present paper is on \emph{stochastic}
 epidemic models.

In the present paper we will study a fairly simple class of stochastic
epidemic models in a closed community, and present properties of the model. These are both small
population properties, and approximations assuming a large community: early stage behaviour of the epidemic, final epidemic size
distribution and the duration of the epidemic. The main large-population approximation results can be
summarized as follows. Assuming a large
population, the early stages of the epidemic can be approximated by a
branching process, where ''giving birth'' corresponds to ''infecting
someone''. If the branching process/epidemic is super-critical it
is possible that a large epidemic outbreak occurs (corresponding to the
branching process growing beyond all limits). If this happens, a balance
equation determines the final number of infected added with some
Gaussian fluctuation of smaller order. As regards to the duration of
a major outbreak the whole
outbreak is divided into three sections: the
beginning (up to when a small fraction have been infected), the main
part (in which nearly all infections take place), and the end (when the
last small fraction of people get infected), and these parts last for durations of
order $\log n$, 1 and $\log n$ respectively, thus making the total
duration of order $\log n$.

We will also describe
how the models can be applied to answer epidemiological questions, for
example how to estimate important
epidemiological parameters from outbreak data and how to study effects
of interventions such as vaccination. We then describe many
important extensions of stochastic epidemic models aiming at making
them more realistic and give some key references. The paper is however
not claiming to be a complete reference guide to all important
contributions in stochastic epidemic models.

In Section \ref{why} we first define the deterministic general
epidemic model and derive some properties of it, then describe some
cases where a deterministic model is insufficient, and end by
defining what we call the standard stochastic SIR-epidemic model. In
Section \ref{mod-prop} we present properties of the model, both
exact for a small population, and approximations relying on a
large community. In Section \ref{appl} we describe how the models can
be used to answers epidemiological questions, and in Section \ref{ext}
we describe a number of model generalizations and also a model for an
endemic infectious disease.

\section{\emph{Stochastic} epidemic models -- why?}\label{why}

\subsection{Deterministic epidemic models}\label{det-mod}

One simple model, the
deterministic general epidemic model (e.g.\ Bailey, 1975, Ch.\ 6.2), can be defined by two
differential equations. It is assumed that at any time point an
individual is either susceptible (s), infected and infectious (i) or
recovered and immune (r). Such individuals are from now one called
susceptibles, infectives and recovered respectively.

The model makes the following assumptions: only susceptible
individuals can get infected and, after having been infectious for
some time, an
individual recovers and becomes completely immune for the remainder of
the study period. Finally, we assume there are no births, deaths,
immigration or emigration during the study period; the community is
said to be \emph{closed}. A consequence of the assumptions is that
individuals can only make two moves: from S to I and from I to
R. For this reason the model is said to be an SIR-epidemic
model. Models having no immunity (individuals that recover become
susceptible immediately) are called SIS-models, models having a
latent state when infected, before becoming infectious, are often
called SEIR (``E'' for exposed but not infectious), models where
immunity wanes after some time are called SIRS-models, and so
forth. Models that allow for births/deaths/immigration/emigration
are referred to as having \emph{demography} or having a \emph{dynamic}
community. The focus in this paper is on SIR models in a closed
community; see however Section \ref{endemicity} for a model allowing births and deaths, and Section \ref{other-ext} for a discussion about latency periods.

Let $s(t)$, $i(t)$ and $r(t)$,
respectively denote the community \emph{fractions} of susceptibles, infectives
and recovered. Since these are fractions and the community is closed we assume that
$s(t)+i(t)+r(t)=1$ for all $t\ge 0$. From the assumptions mentioned
above, together with the assumption of the community being homogeneous
and people mixing homogeneously, the deterministic general epidemic
model is defined by the following set of differential equations:
\begin{align}
s'(t)&= -\lambda s(t)i(t),\nonumber\\
i'(t)&= \lambda s(t)i(t)-\gamma i(t),\label{def-det-gen}\\
r'(t)&= \gamma i(t).\nonumber
\end{align}
These differential equations, together with the starting configuration $s(0)=1-\varepsilon$,
$i(0)=\varepsilon$ and $r(0)=0$ defines the model.

The initial
fraction infectives $\varepsilon >0$ is often assumed to be small as indicated by the
notation $\varepsilon$, it must however be positive -- otherwise all differential equations are constant and equal to 0. The reason for
assuming that $r(0)=0$ is that initially immune individuals play no
part in the dynamics so, up to a normalizing constant, initially immune
individuals may simply be ignored. Some authors choose to let
$s(0)=1$, $i(0)=\varepsilon$ and $r(0)=0$. The use of ``fraction'' is
then somewhat misleading, but $s(t)$ has the interpretation of being
the fraction still susceptible among the initially
susceptibles  at $t$. When $\varepsilon$ is small (as we will nearly always
assume) there is hardly any difference between the two parametrisations.

The term $\lambda s(t)i(t)$ in Equation (\ref{def-det-gen}) comes
from the fact that susceptibles must have contact with infectives in
order to get infected, so the assumption about uniform mixing
(mass-action) implies that infections occur at a rate proportional to
$s(t)i(t)$. This term is non-linear which makes the solution of the
system of differential equations non-trivial. Finally, it is worth
pointing out that since $s(t)+i(t)+r(t)=1$ it is actually enough to
keep track of two of the quantities.

\begin{figure}[!h] \begin{center} \bf
\includegraphics[height=6cm, width=12cm]{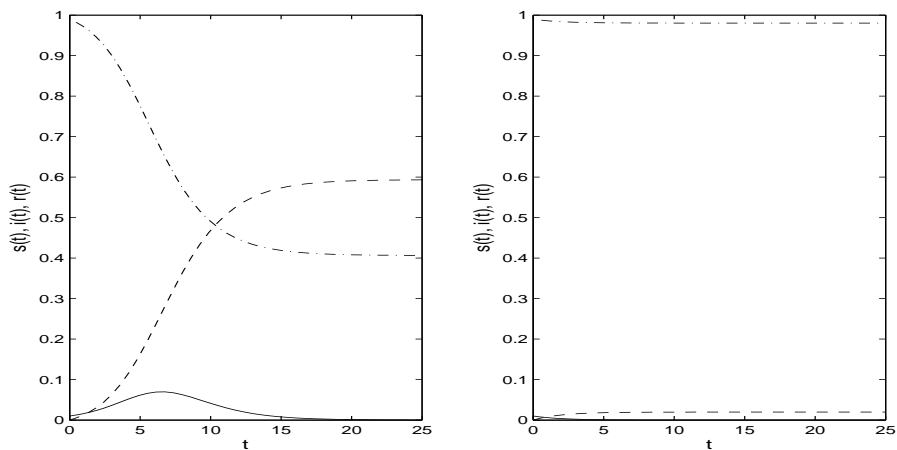}
\caption{\rm Solution of the differential system defined in
  (\ref{def-det-gen}), $s(t)$: $\cdot -$, $i(t)$: ---, $r(t)$: $ -\ -$. Both figures have initial configuration $s(0)=0.99,\
  i(0)=0.01$. To the left is the case with $\lambda=1.5$ and
  $\gamma=1$ (so $R_0=1.5$), and to the right is the case $\lambda=0.5$ and
  $\gamma=1$ (so $R_0=0.5$).}
 \label{fig:diff_gem}\end{center}
\end{figure}

By studying the differential equations it is straightforward to show
that $s(t)$ is monotonically decreasing down to $s(\infty)$ say, and
$r(t)$ is monotonically increasing up to $r(\infty)$. The differential
equation for $i(t)$ can be written as $i'(t)=i(t)(\lambda
s(t)-\gamma)$. So, if $\lambda s(0)>\gamma$, then $i(t)$ initially
increases, but eventually, when $s(t)$ has decreased enough, $i(t)$
starts decreasing. If on the other hand $\lambda s(0)<\gamma$, then
$i(t)$ decreases already from the start with the effect that little
will happen as $t$ tends to infinity (in both cases it can be shown
that $i(\infty)=0$). This dichotomy is illustrated in Figure \ref{fig:diff_gem} where
we have plotted $(s(t), i(t), r(t))$. To the left this is done for the case $\lambda=1.5$, and
$\gamma=1$, and to the right for $\lambda=0.5$ and $\gamma=1$, both
having initial configuration $(s(0),i(0),r(0))=(0.99, 0.01, 0)$. In
the left a substantial fraction (58.3\%) eventually get infected, whereas in
the right figure this fraction is negligible (only an additional 0.9\%
get infected), we say that a major
outbreak has occurred in the first case and a minor outbreak occurred
in the latter case. Since $i(0)$ is assumed small (and $s(0)$ being
close to 1), the critical value separating the two very different
scenarios is $R_0:=\lambda/\gamma =1$.

\begin{figure}[!h] \begin{center} \bf
\includegraphics[height=7cm, width=12cm]{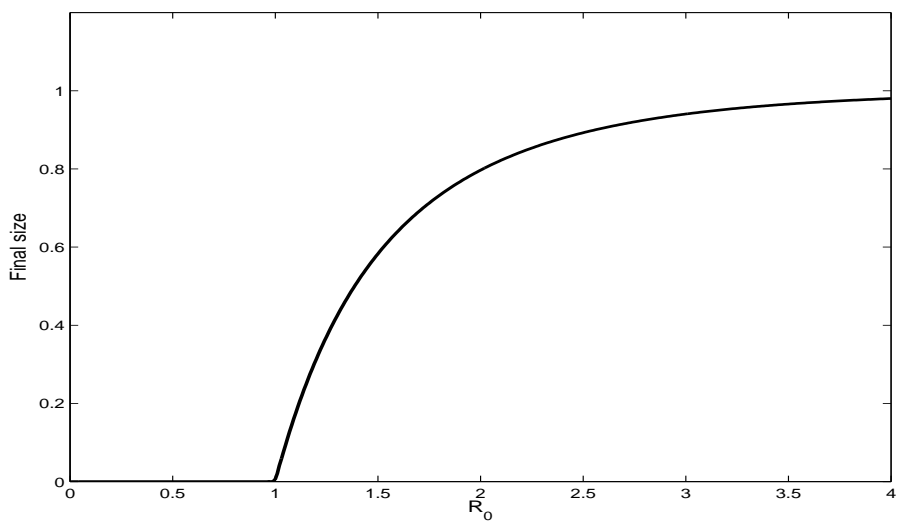}
\caption{\rm The final size of the epidemic as a function of $R_0$ for
the deterministic epidemic model. The initial fraction of infectives
is approximated to equal $i_0=0$.}
 \label{fig:fin_size}\end{center}
\end{figure}

The ratio $R_0=\lambda/\gamma$ is hence of fundamental importance and
can be interpreted as the average number of new infections caused by
an infectious individual before recovering. The ratio is often
referred to as the \emph{basic reproduction number} (a term with its
origin in demography -- the average number of individuals that one
individual reproduces) and denoted by $R_0$:
\begin{equation}
R_0=\frac{\lambda}{\gamma}.\label{R_0}
\end{equation}
When $R_0>1$ the epidemic takes off and when $R_0<1$ there is no (big)
epidemic. The differential equations (\ref{def-det-gen}) can also be
used to obtain a balance equation for the final state $(s(\infty),
0,r(\infty))$. By dividing the first equation by the last we get
$ds/dr=-R_0s$, which implies that $s(t)=s(0)e^{-R_0r(t)}$. The fact
that $i(\infty)=0$ implies that $s(\infty)=1-r(\infty)$; at the end of
the epidemic there are no infectives, only susceptibles and recovered (immune). From this we get a balance equation
determining the fraction $z =r(\infty)$ that at the end of the epidemic
were infected:
\begin{equation}
1-z =(1-\varepsilon)e^{-R_0z }.\label{det-fin-size}
\end{equation}
The balance equation can be interpreted as follows: in order \emph{not} to
have been infected (which a fraction $1-z$ satisfy) you must belong to
those not initially infected (the first factor on the right)
\emph{and} you must escape the infection pressure $R_0$ caused by
those $z$ who were infected. In Figure  \ref{fig:fin_size} we have
plotted the solution $z $ as a function of $R_0$, when starting with
 a negligible fraction initially infectives. The threshold value of 1 is clearly seen to
be the value above which a positive fraction gets infected.

The above results will also be useful when considering a related stochastic epidemic model for a large community.

\subsection{When are deterministic models insufficient?}\label{det-ins}
In the previous subsection we analysed the deterministic general
epidemic model and showed that: if $R_0<1$ there will only be a small
outbreak, and if $R_0>1$ there will be a major outbreak infecting a
substantial fraction of the community, and how big fraction is
determined by Equation (\ref{det-fin-size}). The results rely on that
the community is homogeneous and that individuals mix uniformly with
each other.

Even if the assumption of a homogeneous uniformly mixing community are
accepted this model may not be suitable in some cases. For example, if
considering a small community like an epidemic outbreak in day care
center or school it seems reasonable to assume some
uncertainty/randomness in the final number infected. Also, even if
$R_0>1$ and the community is large but the outbreak is initiated by
only one (or a few) initial infectives it should be possible that, by chance, the epidemic never takes
off. These two arguments motivate the definition of a related \emph{stochastic}
epidemic model. Later we will also show two other reasons motivating
the use of stochastic epidemic models: it enables parameter estimates
from disease outbreak data to be equipped with standard errors and,
when studying epidemic diseases, the question of disease extinction is
better suited for stochastic models.

\subsection{A simple stochastic SIR epidemic model}\label{def-mod}

We now define the \emph{standard stochastic SIR epidemic model}. Just
like for the deterministic general epidemic model we assume a closed
homogeneous uniformly mixing community, and let $n$ denote the
size of the community.

Let $S(t)$, $I(t)$ and $R(t)$ respectively denote the number of
susceptibles, infectives and recovered at time $t$, and suppose that
at time $t=0$ these numbers are given by $S(0)=n-m$, $I(t)=m$ and
$R(0)=0$. The dynamics of the model are defined as follows. Infectious
individuals have "close contact" with other individuals randomly in
time at constant rate $\lambda$, and each such contact is with a
randomly selected individual, all contacts of different infectives
being defined to be mutually independent. By "close contact" is meant
a contact close enough to result in infection if the other individual
is susceptible, otherwise the contact has no effect. Any susceptible
that receives such a contact immediately becomes infected and
infectious and starts spreading the disease according to the same
rules. Infected individuals remain infectious for a random time $I$
(the infectious period) after which they stop being infectious,
recover and become immune to the disease. The infectious periods are
defined to be independent and identically distributed (also
independent of the contact processes) having distribution $F_I$ and
mean $E(I)=1/\gamma$ (to agree with the determinist model).

The epidemic starts at time $t=0$. As the epidemic evolves, according
to the rules above, new individuals (may) get infected and eventually
recover, up until the first time $T$ when there are no infectives in
the community. Then no further individuals can get infected implying that
the epidemic stops. The final state of the epidemic is described by
the ultimate number $R(T)$ infected
(recall that that $I(T)=0$, so $S(T)=n-R(T)$ make up the rest of the
community). The final number of infected $R(T)$ will consist of those
$m$ who were initially infected plus those $Z$, say, who were infected
during the outbreak. Later we will study the exact and approximate
distribution of $Z=R(T)-m$.

Two choices of distributions for the infectious periods have (for
mathematical reasons) received special attention in the
literature. The first is where $F_I$ is exponentially distributed with
intensity parameter $\gamma$, which goes under the name the stochastic
general epidemic model (e.g.\ Bailey, 1975, Ch.\ 6.3). Then the model is
Markovian and the Markov process $(S(t), I(t), R(t))$ has
jump-intensities much related to Equation (\ref{def-det-gen}) of the
deterministic general epidemic. The second choice of infectious period
is where $I$ is non-random (and equal to $1/\gamma$). This choice is
called the continuous-time version of the Reed-Frost model. An
equivalent (for the final outcome) version is where it is assumed that
all infections take place exactly at the end of the infectious period,
and this model was initially defined by Reed and Frost 1928 in a
series of lectures (unpublished). The
Reed-Frost model has the mathematical tractability that whether or not
an individual makes contact with two separate individuals are
independent events. This in turn implies that the Reed-Frost model can
be analysed by an Erd\H{o}s-R\'{e}nyi graph (e.g.\ Bollob\'as, 2001) where each pair
of individual has an edge between them independently, with probability
$p=1-e^{-\lambda /n\gamma}$.

\section{Model properties}\label{mod-prop}
We now explain some important properties of the standard stochastic
SIR epidemic model. In Section \ref{exact} we derive some exact
results; the rest of this section is devoted to
approximations  assuming a large
community.

\subsection{Exact distribution}\label{exact}
No matter what infectious period distribution $F_I$ is considered it
is not possible to derive an exact and simple closed form expressions
for the time dynamics of the epidemic. However, it \emph{is} possible
to derive a recursive formula for the final size of the epidemic,
Equation (\ref{fin-size}) below, and this formula, which we now explain,
is based on the fact that in order not get infected an individual must
"escape" infection from all those who did get infected during the outbreak. This has been
done in several ways (e.g.\ Picard and Lef\`{e}vre, 1990), but our outline
follows that of Ball (1986) where more details can be found.

The derivation of the recursive formula for the final size uses two
main ideas: a Wald's identity for the final size and the total
infection pressure, and the interchangeability of individuals making
it possible to express the probability of getting $i$ additional
infections among the initially $m-n$ in terms of getting $i$ infected
in a smaller subset.

We start with the latter result. To this end, fix $n$ and write $\lambda'=\lambda/n$. Let $Z$ denote the final number infected excluding the initial infectives, so the possible values for $Z$ are $0,\dots ,n-m$. Since individuals are interchangeable we can label the individuals according to the order in which they get infected. The initial infectives are labelled $-(m-1),\dots, 0$, then according to time of infection: $1,\dots , Z$, and those who avoid infection according to any order $Z+1, \dots ,n-m$. With this labelling we define the total infection pressure $A$ by
\begin{equation}
A=\lambda'\sum_{i=-(m-1)}^ZI_i \label{tot-inf-pres}
\end{equation}
i.e.\ the infection pressure, exerted on any individual, during the outbreak (sometimes also called the "total cost" of the epidemic).

Now, let $p_i^{(n-m)}=P(Z^{(n-m)}=i)$ denote the probability that exactly $i$ susceptibles get infected during the outbreak, explicitly showing the number of initial susceptibles but suppressing the dependence on the initial number of infectives $m$.
Then, using the interchangeability of individuals and reasoning in terms of subsets among the initially susceptibles, it can be shown (Ball, 1986) that for any $i\le k\le n-m$, it holds that
\begin{equation}
\frac{p_i^{(n-m)}}{\binom{n-m}{i}} = \frac{p_i^{(k)}}{\binom{k}{i}} E\left(e^{-(n-m-k)A^{(k)}}|Z^{(k)}=i\right);\label{subset-prob}
\end{equation}
the probability that $i$ get infected among the initially $m-n$ initially susceptibles equals the product of the probability of having $i$ infected in a smaller subset of size $k$ ($k\le m-n$) multiplied by the probability that no one in the remaining set gets infected conditioning on the first event. The notation $A^{(k)}$ and $Z^{(k)}$ are hence the total pressure and final size starting with $k$ susceptibles (and $m$ infectives).

We now show a Wald's identity for $Z^{(k)}$ and $A^{(k)}$. Let
$\phi(\theta)=E(e^{-\theta I})$ denote the Laplace transform of the infectious period $I$. We then have the following Wald's identity (Ball, 1986):
\begin{equation}
E\left(\frac{e^{-\theta A^{(k)}}}{(\phi (\theta \lambda'))^{m+Z^{(k)}}}\right)
= 1,\qquad \theta \ge 0.\label{Wald}
\end{equation}
The following steps proves the result
\begin{align*}
(\phi (\theta \lambda'))^{m+k}&= E\left[ \exp\left(-\theta \lambda'\sum_{i=-(m-1)}^{k}I_i \right)\right]
\\
&= E\left[ \exp\left(-\theta \left(A^{(k)}+\lambda'\sum_{i=Z^{(k)}+1}^{k}I_i\right) \right)\right]
\\
&= E\left[ e^{-\theta A^{(k)}}(\phi (\theta \lambda')^{k-Z^{(k)}} \right],
\end{align*}
where the last identity follows because the $k-Z^{(k)}$ infectious periods
$I_{Z^{(k)}+1},\dots I_{k}$, are mutually independent and
also independent of $A^{(k)}$ (which only depends on the first $Z^{(k)}$ infectious periods and the contact processes of these individuals).
Dividing both sides by $(\phi (\theta \lambda'))^{m+k}$ gives the desired result.

If we use Wald's identity with $\theta =n-m-k$ and condition on the value of $Z^{(k)}$ we get
\begin{equation}
\sum_{i=0}^k\frac{E\left(e^{-(n-m-k)A^{(k)}}|Z^{(k)}=i\right)}{(\phi
  ((n-m-k) \lambda'))^{m+i}}p_i^{(k)}=1.
\label{wald-cond}
\end{equation}
Using Equation (\ref{subset-prob}) in the equation above we
get
\[
\sum_{i=0}^k \frac{\binom{k}{i}p_i^{(n-m)}}{\binom{n-m}{i}(\phi
  ((n-m-k) \lambda'))^{m+i}} =1.
\]
Simplifying the equation, returning to $\lambda=\lambda' n$ and
putting $p_k^{(m-n)}$ on one side, we obtain the recursive formula for the final size distribution $p_k^{(n-m)}, k=0,\dots ,n-m$:

\textbf{Exact final size distribution}:
\begin{equation}
P_k^{(n-m)} = \binom{n-m}{k}[\phi ((n-m-k) \lambda/n)]^{m+k}
-\sum_{i=0}^{k-1} \binom{n-m-i}{k-i}[\phi
  ((n-m-k) \lambda/n)]^{k-i}p_i^{(n-m)}.
  \label{fin-size}
\end{equation}
Note that this is a recursive formula. For example,
solving Equation (\ref{fin-size}) for $k=0$ and then for $k=1$ gives, after some algebra,
\begin{align*}
p_0^{(n-m)}&=(\phi(\lambda))^m,
\\
p_1^{(n-m)}& = n\phi((n-1)\lambda/n)\Big( [\phi((n-1)\lambda/n)]^{m} -
[\phi(\lambda)]^{m}\Big).
\end{align*}
In order to compute
$p_k^{(n-m)}$ using (\ref{fin-size}) it is hence necessary to sequentially compute $p_0^{(n-m)}$ up to
$p_{k-1}^{(n-m)}$. As a consequence the formula is not very
enlightening and it may be numerically very unstable when $k$ (and
hence $n-m\ge k$) is large. Even when it is possible to compute $p_{k}^{(n-m)}$ using
 (\ref{fin-size}) with $n-m$ being large, an
approximative formula can be more informative. For this reason we
devote the rest of this section to the case where $n-m$ is large.

\subsection{Early stage approximation}\label{br-appr}
We start by deriving an approximation for the early stages of an
outbreak. The approximation relies on that the initial number $n-m$ of
susceptibles is large. The key reason for having an approximation
during the early stages of an outbreak when there are many initial
susceptibles is that it is then very unlikely that any of the first
number of infectious contacts happen to be with the same susceptible
individual. Conversely, it is very \emph{likely} that all of the first
set of infectious contacts happen with distinct individuals. As a
consequence, the number of individuals that distinct infectives infect
(during the early stages) are independent and identically
distributed random variables. This hence suggests that the epidemic
process may be approximated by a suitable branching process (e.g.\
Haccou et al., 2005). In this subsection we make the approximation more precise and derive results for when a major outbreak (infecting a non-negligible fraction) is possible, and if so, also what the probability of a major outbreak is.

Assume that the initial number of infectives $m\ge 1$ is fixed and
that $n-m$ is large (later we will take limits as $n\to\infty$). One
way of constructing the epidemic is as follows. Label the $n$
individuals $-(m-1),\dots , n-m$, where the first $m$ individuals
refer to the initial infectives and the remaining $n-m$ refer to the
initial susceptibles, but otherwise arbitrary -- note that here the
labelling is not according to order of infection. Let
$I_{-(m-1)},\dots ,I_{n-m}$ be i.i.d.\ with distribution $F_I$ (the
infectious periods), $\chi_{-(m-1)} (\cdot ), \dots, \chi_{n-m} (\cdot
)$ being i.i.d.\ Poisson processes all having constant intensity
$\lambda$ (the contact processes), and let $U_1,\dots ,U_{n-m}$ be
i.i.d.\ uniform random variables on the unit interval (to be used for
determining who is contacted when a close contact occurs). The epidemic process is defined using
these random variables as follows. Start at $t=0$. The contact
processes of the initial infectives $\chi_{-(m-1)} (\cdot ), \dots,
\chi_{0} (\cdot )$ are "activated" and the infectious periods of these
individuals start. Time increases without anything happening until the
first time point $t_1$ when either one of the activated contact
processes has an "event", or one of the infectious periods
$I_{-(m-1)},\dots ,I_{0}$ stops. If the latter happens the
corresponding contact process is deactivated and that persons
infectious period stops and the individual recovers and becomes immune. If the former happens the corresponding
individual has an infectious contact. Which person that is contacted
is determined by $U_1$: the person being contacted has index $\lfloor
nU_1\rfloor -(m-1)$, where $\lfloor \cdot \rfloor$ denotes the integer part, implying that each individual $i$ ($i=-(m-1),\dots , n-m$) is selected with equal probability $1/n$. If the contacted person is still susceptible he/she gets infected and the corresponding infectious period is started and the contact process activated, and if the contacted person has already been infected nothing happens. Time then moves on until either an infectious period is terminated or a contact among the activated contact processes occur. Depending on what happens an infectious period is terminated or a randomly selected individual is contacted (and infected unless it has already been infected). This goes on until the first (random) time point $T$ when there are no activated contact processes and all initiated infectious periods have stopped. Since there are finite number of individuals, all having finite infectious periods with probability 1, $T$ will be finite with probability 1.

It is straightforward to check that this construction gives the desired epidemic model: individuals have infectious periods distributed according to $F_I$ and while infectious each individual has contacts with randomly selected individuals at rate $\lambda$. Another nice feature with this construction is that it is in fact possible to construct a sequences of epidemics, indexed by the initial number of infectives, as well as a homogeneous branching process with life-lengths distributed according to $F_I$ and constant birth rate $\lambda$, using the same set of random variables. This can be used to \emph{couple} the whole sequence of epidemics and the limiting branching process and to show that they agree up to some random point. We refer to Ball (1983) for a more detailed study about this -- here we just give the heuristics.

The branching process is simply constructed without the uniform
numbers, so in the branching process a new individual is "born"
whenever a contact occurs. The same applies to the epidemic (having
$n-m$ initially susceptible) with "born" replaced by "infected" except
when a contact in the epidemic is with an already infected
individual. As a consequence, the epidemic and the branching process
agrees up until the first time point when a contact is with an already
infected in the epidemic, denoted a ghost contact in Mollison (1977). Initially
there are $n-m$ susceptible and $m$ infected, so the probability that
the first contact is \emph{not} a ghost contact equals $(n-m)/n$ which
is close to 1 when $n$ is large. Given this, the second contact is
also a non-ghost-contact with probability $(n-m-1)/n$, and so on; note the
resemblance with the well-known birthday problem. From this it
follows that the branching process and epidemic agrees at least up
until the $k$'th contact (i.e.\ no ghost contact has occurred) with probbility
\begin{equation}
P(\text{no ghost among first $k$ contacts})=\frac{(n-m)_k}{n^k}, \label{exactghost}
\end{equation}
where $r_j:=r(r-1)\dots (r-j+1)$. Recalling that $m$ and $k$ are much smaller than $n$ and using the well-known approximation $1-\varepsilon\approx e^{-\varepsilon}$, we get the following approximation for the probability in (\ref{exactghost}):
\[
P(\text{no ghost among first $k$ contacts}) \approx e^{-\left(\frac{m}{n}+\dots + \frac{m+k-1}{n}\right)} =e^{-k(k-1+m)/2n}.
\]
For large $n$ this probability is close to 1 whenever $k=o(\sqrt{n})$. From this it follows that, with large probability, the epidemic and the branching process agrees at least up until there has been $k$ contacts, where $k$ is small in relation to $\sqrt{n}$.

The above heuristics (made precise in Ball, 1983) motivate that we can
approximate the epidemic with the branching process up until $k$
births have occurred. The advantage with this approximation is that
branching processes are well-studied (e.g.\ Haccou et al., 2005). For
instance, our branching process (with life-length distribution $I\sim
F_I$ and constant birth rate $\lambda$) has mean offspring
distribution $\lambda E(I)=\lambda/\gamma$, a quantity previously
denoted by $R_0$ ($=\lambda/\gamma$). The branching process is
sub-critical if $R_0<1$, critical if $R_0=1$ and supercritical if
$R_0>1$. The total progeny $Z_\infty$, the number of individuals ever
born in the branching process, is known to be finite (the branching
process goes extinct) with probability 1 if $R_0\le 1$. And,
if $R_0>1$ then $Z_\infty$ has a finite part, and is infinite (grows
beyond all limits) with a strictly positive probability $\rho$ that
can be determined. In either case, the distribution of the finite part
of $Z_\infty$ has a well-defined distribution depending on $m$,
$\lambda$ and $F_I$, cf.\ Jagers (1975, Ch 2.11).

Due to the approximation outlined above, the branching process and the epidemic coincide on the part of the sample space where the branching process goes extinct (whether sub-critical or not), with arbitrary large probability if $n$ is large enough. As a consequence, we may approximate $P(Z_n=j)$ with $P(Z_\infty =j)$ for small and moderate $j$, a distribution which is quite complicated except for some special cases. To compute an expression for $\rho =P(Z_\infty =\infty)$ is however easier. We compute the opposite, i.e.\ the probability $1-\rho$ that the branching process, starting with $m$ individuals, goes extinct. For this to happen all $m$ initial lineages must go extinct. So if $q$ denotes the probability that one lineage goes extinct it follows that $1-\rho =q ^m$. In order to obtain an expression for $q$ we condition on the number of births $X$ the initial individual gives birth to before dying. Given that $X=j$ all these $j$ offspring must have lineages that go extinct, an event that happens with probability $q^j$. It follows that $q$ must satisfy $q=E(q^X)$. The number $X$ that an individual gives birth to during the life-length $I\sim F_I$ of course depends on $I$. Given this duration: $I=y$, the number of births is Poisson distributed with mean parameter $\lambda y$ since contacts occur according to a Poisson process having rate $\lambda$. So, if we denote the Laplace transform of the infectious period $I$ by $\phi (s)=E(e^{-sI})$, then the following relation must hence hold for $q$:
\begin{equation}
q=E(q^X) =E(E(q^{X}|I)) = E(e^{-\lambda I(1-q)})=\phi (\lambda(1-q)),\label{q-def}
\end{equation}
where the third equality uses the Laplace transform of a Poisson variable.
It is in fact known (e.g.\ Jagers, 1975, p 9) that $q$ is the smallest solution to this equation.

\textbf{Early stage approximation}.
To summarize, when $n$ is large, the initial phase of the epidemic can be approximated by a homogeneous branching process with birth rate $\lambda$ and life-distribution $I\sim F_I$ having Laplace transform $\phi(s)$. If $R_0:=\lambda E(I)=\lambda /\gamma \le 1$ the final size of the epidemic is bounded in probability, whereas if $R_0>1$ it is not. The approximation further tells us that when $R_0>1$ the epidemic will be "minor" (bounded) with a probability $q^m$ and "major" (unbounded) with probability $\rho=1-q^m$, where $q$ is the smallest solution to (\ref{q-def}). The distribution of outbreak sizes in the "minor" case can also be determined from branching process theory.

We end this section with an example illustrating our results. In the section to follow we return to the situation where the branching process grows beyond all limits, corresponding to a major outbreak in the epidemic.

\textbf{Example}.  Suppose the epidemic is initiated with $m=1$
individual and that there are $n-1=999$ initial susceptibles. Further
we treat the case where the infectious period is $F_I=Exp (\gamma)$,
the exponential distribution with rate parameter $\gamma$ (i.e.\ the general stochastic epidemic). When $I\sim {\rm Exp}(\gamma)$ then
$\phi (s)=E(e^{-sI})=\gamma/(\gamma +s)$. This implies that $q$ is the
smallest solution to $q=\gamma/(\gamma +\lambda (1-q))$. Solving this
quadratic equation shows that $q=\gamma/\lambda =1/R_0$ if $R_0>1$ and
otherwise $q=1$. This means that there will be a minor epidemic
outbreak with probability $1/R_0$, and a major outbreak with the
remaining probability $1-1/R_0$. Further, it can be shown that for this particular infectious period/life-length distribution (computations omitted) the outbreak probabilities are well approximated by the corresponding total progeny distribution
\[
P(Z_n=j)\approx \binom{2j}{j}\frac{1}{(j+1)!}\left(\frac{1}{1+R_0}\right)^{j+1} \left(\frac{R_0}{1+R_0}\right)^{j},\quad j=0,1,\dots .
\]
Note that for fixed $\lambda$ and $\gamma$ (and hence
$R_0=\lambda/\gamma$) both $q$ and $P(Z_n=j)$ can be computed
numerically. However, the approximation for $P(Z_n=j)$ relies on that
not too many individuals have been infected -- then the branching
process approximation breaks down -- so it should only be used for $j$
up to around 20 ($\sqrt{n}\approx 31.6$).

\subsection{Final size approximation}\label{sec-stoch-fin-size}

In the previous subsection it was shown that when $n$ is large the epidemic is well-approximated by a branching process up until approximately $\sqrt{n}$ individuals have been born. If the branching process goes extinct this will never happen when $n$ is large enough, but when the branching process grows beyond all limits (hence implicitly assuming that $R_0>1$) the approximation breaks down. The question is of course what happens with the epidemic in this case; something which we now briefly outline. The outline is based on results in Scalia Tomba (1985, 1990) which uses the elegant Sellke construction (Sellke, 1983).

 When we are only interested in the final number infected we may introduce a different time scale as follows. Each initially susceptible individual is given a so-called \emph{resistance} having ${\rm Exp}(1)$ distribution. In the first time step we let each initially infective $i$ "throw out" its infection pressure $\lambda I_i$ (defined as the contact rate multiplied by the length of the infectious period) uniformly in the community (so $\lambda I_i/n$ on each individual). Those individuals with resistances smaller than this infection pressure then become infected and in turn throw out their infection pressure uniformly, thus increasing the accumulated infection pressure. This procedure goes on until the first time step when there are no new infections; then the epidemic stops. It can be shown that this gives the correct final size distribution -- the only difference with this representation lies in the order and time at which individuals get infected. The advantage with the construction (the Sellke construction) is that, at each time step we add a random number of i.i.d.\ random variables (the infection pressures), and this is done until the first time point the accumulated sum no longer exceeds another sum of random variables (the resistances). Scalia Tomba (1985, 1990) uses this, together with an embedding argument, to show that the final size distribution coincides with that of the crossing time of a stochastic process with a straight line. Because the random process is made up of i.i.d.\ contributions it obeys a law of large number and central limit theorem.

The above description is perhaps not very enlightening, but to show the complete result is quite technical. Something which is easier to explain is an argument for the expression of the limiting final fraction infected $z$ in case there is a major outbreak. We do this by deriving a balance equation for $z$. Neglecting the difference between $n$ and $n-m$ (remember that $n$ is large and $m$ fixed), we have that $nz$ is the final \emph{number} infected. Further, the probability of escaping infection from one infected individual $i$ (with infectious period $I_i$) equals $E(e^{-\lambda I_i/n})$. We then have the following set of approximations
\begin{align*}
1-z& =\text{fraction not getting infected}\\
&\approx\text{probability of not getting infected}\\
&\approx\text{probability of escaping infection from all $nz$ getting infected}\\
&\approx E(e^{-\lambda I_1/n})\cdot \dots \cdot E(e^{-\lambda I_{nz}/n})\\
&=E\left(e^{-\frac{\lambda}{n}\sum_{i=1}^{zn}I_i} \right)\\
&\approx e^{-\lambda z E(I)} =e^{-R_0z},
\end{align*}
where the last approximation is simply the law of large numbers. The limiting fraction infected in case there is a major outbreak should hence be a solution to the equation
\begin{equation}
1-z=e^{-R_0z}.\label{fin-size-eq}
 \end{equation}
Note that this is the same equation as for the deterministic model
(\ref{det-fin-size}), except that we now assume a negligible fraction
initially infected (i.e.\ $\varepsilon =0$).
 It is not hard to show that this equation always has the solution
 $z=0$ (corresponding to a minor outbreak) and, when $R_0>1$, there is
 another unique solution $z^*$ between 0 and 1 (corresponding to a
 major outbreak).  In Figure \ref{fig:fin_size} of Section \ref{def-mod}, the largest solution $z^*$ is plotted as
 a function of $R_0$.

The above heuristics indicates that the final fraction infected $Z/n$ will lie
close to either 0 or else, if $R_0>1$, close to $z^*$. This can in
fact be shown rigorously. Moreover, a central limit theorem can be
shown for the case where there is a major outbreak. The following
theorem summarizes the result for both a minor and a major outbreak
(see Andersson and Britton, 2000a, Ch 4, and references therein, e.g.\
von Bahr and Martin-L\"of, 1980).

\textbf{Threshold theorem for the final size of the epidemic}: 
Consider the standard epidemic model with $m$ (fixed) initial
infectives and $n-m$ initial susceptibles. If $R_0\le 1$, then $\bar
Z_n:=Z_n/n\to 0$ in probability as $n\to\infty$. On the other hand, if
$R_0>1$, then $\bar Z_n\to \zeta$ in distribution, where $\zeta$ is a
two point distribution defined by $P(\zeta=0)=q^m$ and
$P(\zeta=z^*)=1-q^m$, where $q$ was defined in the previous subsection
and $z^*$ is the unique strictly positive solution to
(\ref{fin-size-eq}). Finally, on the part of the sample space where
$\bar Z_n\to z^*$ we have that
\begin{equation}
\sqrt{n}\left( \bar Z_n -z^*\right) \to N\left( 0,\
  \frac{z^*(1-z^*)\left(1+r^2(1-z^*)R_0^2\right)}{(1-(1-z^*)R_0)^2}\right) ,\label{clt}
\end{equation}
where the second parameter in the normal distribution denotes the
variance, in which $r^2=V(I)/(E(I))^2$ is the squared coefficient of
variation of the infectious period.

To illustrate our results we have simulated an epidemic in a community
of $n=1000$ individuals, starting with $m=1$  initially infective, and having
$\lambda=1.5$ and $I\sim {\rm Exp} (1)$. Before looking at the simulations
we compute the theoretical values. For our model and parameter choices
we have $R_0=\lambda E(I)=1.5$ and $r^2=V(I)/(E(I))^2=1/1=1$. Using
results from the previous subsection we conclude that the probability
of a minor outbreak equals  $q^1=1/R_0\approx
0.667$, and using (\ref{fin-size-eq}) we find that $z^*\approx 0.583$. The standard deviation of $Z_n$ equals
$\sqrt{n}\sqrt{z^*(1-z^*)\left(1+r^2(1-z^*)R_0\right)}/(1-(1-z^*)R_0)
\approx 53.1$.

\begin{figure}[!h] \begin{center} \bf
\includegraphics[height=!, width=10cm]{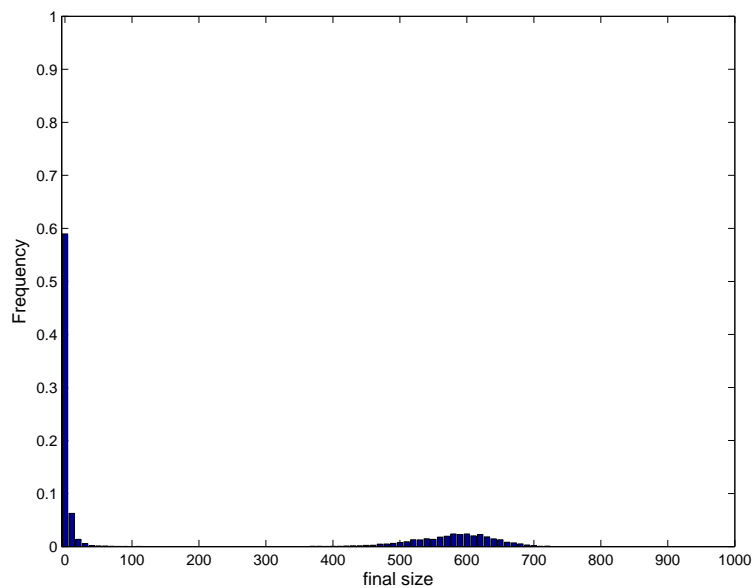}
\caption{\rm Empirical distribution of final size from 10~000
  simulations of the general stochastic epidemic with $R_0=1.5$ in a
  community with 1000 individuals and one initial infective.}
 \label{fig:fin-size-histo}\end{center}
\end{figure}

To see how our limiting results apply to this particular finite case
we
have simulated 10~000 such epidemics. In Figure
\ref{fig:fin-size-histo} we show the resulting histogram for the final
outbreak sizes. The conclusion is that the limiting results apply quite
well with the simulations. First,
the proportion of simulations resulting in minor outbreaks (as seen from the figure there is a clear distinction between minor and major outbreaks) equals 0.678, which is to be compared with the theoretical
value of $2/3\approx 0.667$. Further, by simple examination it looks
as if the remaining part (the major outbreaks) have a normal
distribution with a mean close to
$nz^*\approx 583$ as suggested by the limiting result. Finally, the
standard deviation is of course harder to ''see'', but nearly all
observations seem to lie within 100 from the mean, which agrees well
with the theoretical value of 53 for the standard deviation.

\subsection{Duration of epidemic}\label{duration}

In the previous sections we have studied the questions of main interest: can an outbreak occur, and if so, with what probability and how large will these major outbreaks be. Another questions of interest could be to understand how long it will take for the epidemic to peak and eventually to die out. Below we sketch some results in his direction. For details we refer to Barbour (1975).

If there is only a minor outbreak (which happens with certainty if $R_0\le 1$, but also with positive probability $q^m$ if $R_0>1$), then the disease will disappear after a short time. We hence focus on the case where we have a major outbreak (hence assuming $R_0>1$) and study how the time to extinction $T=T_n$ depends on the population size $n$.

It was seen that the initial stages of the epidemic could be approximated by a supercritical branching process up until approximately $\sqrt{n}$ individuals have been infected. Since a branching process has exponential growth, this will take a time of order $\log n$. Once a large number of individuals have been infected, the epidemic process may be approximated by the deterministic counter part defined in (\ref{def-det-gen}). The time for this process, starting in an arbitrary small fraction initially infectives $i_0$ to first grow and the decrease down to some small value $i(t)=\varepsilon$ has a duration which does not grow with $n$. The last part of the epidemic is where close to a fraction $z^*$ have already been infected. The epidemic then behaves like a branching process where individual infect on average $R_0z^*$ individuals. It can be shown that this number is smaller than 1 when $R_0>1$ as we have assumed. It follows that the end of the epidemic can be approximated by a subcritical branching process starting with $n\varepsilon$ infectives. The duration for such a process to go extinct is also of order $\log n$.

To summarize, we have made it plausible that the duration of the epidemic $T_n$, assuming a major outbreak, has a distribution of the form
\begin{equation}
T_n=c_1\log n + c_2 + c_3\log n + Z,\label{T_n}
\end{equation}
 where $c_1,\ c_2,\ c_3$ are constants and $Z$ is a random variable, all depending on the model parameters. To show this result is quite technical, see Barbour (1975).

 \section{Applications}\label{appl}
 There are many applications to modelling of infectious disease spread. Some diseases that have received much modelling attention over the last decades are for example HIV, Small-pox (the threat for terrorist attacks), foot and mouth disease, SARS and most recently the new (H1N1) influenza. Below we describe some specific methodological questions that are often addressed.

\subsection{Vaccination and other interventions}\label{vacc}

One reason for modelling infectious disease spread is to understand how an outbreak can be prevented. This can be achieved in several different ways. When it comes to new emerging severe infections like SARS, drastic measures like isolation, closing of schools and travel restrictions are often put into place. All these measures aim at reducing contact rates, i.e.\ to reduce $R_0$. The effect of a specific preventive measure depends on the particular disease and also on the community under consideration.

A somewhat different preventive measure is vaccination. This does not change $R_0$, but instead it reduces the pool of susceptibles by making individuals immune. We now describe how to model this and study its effect on the outbreak.

Suppose a vaccine is available prior to the arrival of the disease, and assume a fraction $v$ are vaccinated. Suppose further that the vaccine is perfect in the sense that all vaccinated individuals get completely immune (see Section \ref{other-ext} for extensions).

This will have one effect on the model: the number of initially susceptibles is reduced from $n-m$ to $n(1-v)-m=:n'-m$. 
Other than that the model is the same. However, the contact rate with a given individual still equals $\lambda/n$, so if we want to use $n'$ instead we get $\lambda/n=\lambda (1-v)/n'=:\lambda'/n'$.
We can hence compute the exact final size distribution using results in Section \ref{exact}, for example Equation (\ref{fin-size}) with $n'=n(1-v)$ instead of $n$ and $\lambda'=\lambda (1-v)$ instead of $\lambda$. Similarly, the early stages of the epidemic can be approximated by a branching process with birth rate $\lambda'$ and mean life-length $1/\gamma$. It follows that the reproduction number after a fraction $v$ has been vaccinated, denoted $R_v$, satisfies
\begin{equation}
R_v=\lambda '/\gamma =\lambda (1-v)/\gamma =(1-v)R_0.\label{R_v}
\end{equation}

From the results of Section \ref{br-appr} we hence conclude that if $R_v\le 1$ there will be no major outbreak, whereas if $R_v>1$ there will be a major outbreak with probability $1-(q')^m$, where $q'$ is the smallest solution to Equation (\ref{q-def}), but with $\lambda'=\lambda (1-v)$ instead of $\lambda$.

From Section \ref{sec-stoch-fin-size} we also conclude that if $R_v>1$, the major outbreak will be approximately $n'z'^*$, where $z'^*$ is the unique positive solution to Equation (\ref{fin-size-eq}) with $R_v$ replacing $R_0$. The central limit theorem also applies, so in case there is a major outbreak the total number of infected is normally distributed as stated in (\ref{clt}), but wit $n'$ and $\lambda'$ and $R_v$ instead of $n$, $\lambda$ and $R_0$.

The most important of these results from an applied point of view is (\ref{R_v}), i.e.\ that the reproduction number after having vaccinated a fraction $v$ in the community satisfies $R_v=R_0(1-v)$, and the fact that a major outbreak is impossible if $R_v\le 1$. In terms of $v$ this is equivalent to $v\ge 1-1/R_0$. The critical vaccination coverage, denoted $v_C$ and defined as the fraction necessary to vaccinate in order to surely prevent a major outbreak, hence satisfies
\begin{equation}
v_C=1-\frac{1}{R_0}.\label{v_C}
\end{equation}
For the numerical example given above, with $R_0=1.5$ it hence follows that $v_C=1-1/1.5\approx 0.33$. This means that it is enough to vaccinate 33\% of the community to prevent a major outbreak. By vaccinating only 33\% the whole community is hence protected, a state denoted \emph{herd immunity}.

\subsection{Estimation}\label{inf}
So far we have been interested in (stochastic) modelling, i.e.\ given a model and its parameters we have studied properties of the epidemic. It was shown that the most important parameter is the basic reproduction number $R_0$, defined as the average number of individuals a (typical) infective infects during the early stages of the outbreak, where we by ''early stages'' mean that few individuals have been infected so the vast majority of the community is still susceptible. When aiming at preventing an outbreak, another important quantity is the critical vaccination coverage $v_C$ defined in (\ref{v_C}).

In a particular situation it is hence important to know what these parameters are. We now describe how to estimate $R_0$ and $v_C$ from observing one epidemic outbreak. One advantage of stochastic modelling is that it not only enables point estimates, but also standard errors as is now illustrated.

Suppose a major epidemic outbreak occurs in a community of $n$ individuals, $n-m$ being susceptible and a small number $m$ being externally infected. Let $z_n$ denote the number of individuals that were infected during the outbreak (excluding the $m$ index cases), and let $\bar z_n=z_n/n$ denote the corresponding fraction. We want to estimate $R_0$ (and $v_C$) based on this information. From (\ref{fin-size-eq}) we know that $R_0$ satisfies
$$
R_0=\frac{-\log (1-z^*)}{z^*}.
$$
We also know that, in case of a major outbreak, $\bar Z_n$ is asymptotically normally distributed around $z^*$, with standard deviation
$$
\frac{\sqrt{z^*(1-z^*)\left(1+r^2(1-z^*)R_0\right)}}{\sqrt{n}(1-(1-z^*)R_0)}.
$$
From this it follows that a consistent and asymptotically normally distributed estimator of $R_0$ is given by
\begin{equation}
\hat R_0= \frac{-\log (1-\bar z_n)}{\bar z_n}.\label{R_0-est}
\end{equation}
The standard deviation of the estimator can be obtained using the
delta method (e.g.\ Cox, 1998). Let $f(z)=-\log (1-z)/z$, then
$f'(z^*)=1/(z^*(1-z^*))-R_0/z^*$. It follows that the asymptotic
variance of $\hat R_0$ hence equals
\begin{align*}
V(\hat R_0)& \approx (f'(z^*))^2V(\bar Z_n) \approx \left(\frac{1}{z^*(1-z^*)} - \frac{R_0}{z^*}\right)^2 \frac{z^*(1-z^*)\left(1+r^2(1-z^*)R_0^2\right)}{n(1-(1-z^*)R_0)^2 }\\
& =\frac{ 1+r^2(1-z^*)R_0^2}{nz^*(1-z^*)}.
\end{align*}
A standard error for $\hat R_0$ is obtained by taking the square root of this and replacing $R_0$ by $\hat R_0$ and $z^*$ by $\bar z_n$:
\begin{align}
s.e.(\hat R_0)= \sqrt{\frac{ 1+r^2(1-\bar z_n)\hat R_0^2}{n\bar z_n (1-\bar z_n)}}. \label{se-R_0}
\end{align}
The standard error above still contains one unknown quantity: $r^2=V(I)/(E(I))^2$, the squared coefficient of variation of the infectious period. From final size data it is impossible to infer anything about the distribution of the infectious period. The only way to proceed, unless some prior information is available, is to replace it by some conservative upper bound, for example $r^2=1$.

Estimation of the critical vaccination $v_C$ is also straightforward. Since $v_C$ is defined by (\ref{v_C}) the natural estimator for $v_C$ is
\begin{align}
\hat v_C=1-\frac{1}{\hat R_0} = 1- \frac{\bar z_n}{-\log (1-\bar z_n)}.\label{v_C-est}
\end{align}
Just like $\hat R_0$ the estimator $\hat v_C$ is consistent and asymptotically normally distributed around the true value $v_C$. The asymptotic variance can also be obtained using the delta method. We have $g(x)=1-1/x$, which hence satisfies $g'(R_0)=1/R_0^2$. The asymptotic variance of $\hat v_C$ then equals $V(\hat R_0)/R_0^4$. A standard error for $\hat v_C$ is given by
\begin{align}
s.e.(\hat v_C)= \sqrt{\frac{ 1+r^2(1-\bar z_n)\hat R_0^2}{n\hat R_0^4\bar z_n (1-\bar z_n)}}. \label{se-v_C}
\end{align}

As a numerical example we assume that $z_n=583$ individuals in a community of $n=1000$ were infected during an outbreak. Using (\ref{R_0-est}) and (\ref{se-R_0}) we get $\hat R_0 =1.50$ with standard error $s.e.(\hat R_0)=0.09$ using the conservative upper bound $r^2=1$. As for the critical vaccination coverage the corresponding estimate and standard error are given by $\hat v_C=0.333$ and $s.e.(\hat v_C)= 0.04$. Both estimators are asymptotically normally distributed.

In the inference procedure presented above we only made use of the final number infected. If the epidemic process $(S(t), I(t), R(t))$ is observed continuously we have more detailed information and we should hence be able to make more precise inference. This is true although the gain in precision is moderate. We omit this type of inference (involving martingales using counting process theory) and refer to Andersson and Britton (2000a, Ch 9).

\section{Extensions}\label{ext}
The standard stochastic epidemic model of Section \ref{def-mod} was defined as a stochastic version of the deterministic general epidemic model of Section \ref{det-mod}. The population is finite, and infectious individuals make contacts randomly in time according to a Poisson process with intensity $\lambda$, each time the person being contacted is chosen randomly, and the length of the infectious period $I$ is a random variable with distribution $I$. Despite these more realistic model features, the model still contains several simplifying assumptions. In the present section we touch upon a number of extensions that have been analysed in the literature. Most of these extensions do not alter the qualitative behaviour of the spread of infection, but indeed quantitative properties.

\subsection{Individual heterogeneity: multitype models}\label{multitype}
One assumption in the standard stochastic epidemic model is that all
individuals are similar, except the possibility that they have
different length of the infectious period. In reality most populations
are heterogeneous, for example with respect to susceptibility, the
degree of social activity and/or how infectious they become if
infected. The heterogeneities may be unknown to some extent, but quite
often it is possible to group individuals into different \emph{types}
of individual, where individuals of the same type have (nearly)
identical behaviour. This separation into different types might for
example be different age groups, gender, previous experience to the
disease (giving partial immunity) etcetera.

In such situations it is common to define a multitype epidemic model as follows. Suppose there are $k$ types of individuals, labelled $1,\dots ,k$, and that the community fractions of the different types equal $\pi_1,\dots ,\pi_k$. Type $i$ individuals have i.i.d.\ infectious periods $I_i$ with distribution $F_i$ having mean $1/\gamma_i$. During the infectious period, an $i$-individual has close contacts with a given $j$-individual at rate $\lambda_{ij}/n$, $i,j=1,\dots ,n$.

If the population, i.e.\ $n$, is large, the epidemic may be approximated by a multitype branching process. The mean offspring matrix has elements $(m_{ij})$, where $m_{ij}=\lambda_{ij}\pi_j/\gamma_i$. The basic reproduction number $R_0$ is the largest eigenvalue to the mean offspring matrix $(m_{ij})$. The case where $m_{ij}=\alpha_i\beta_j\pi_j$, referred to as proportionate, or separable mixing, cf.\ Diekmann and Heesterbeek (2000, Ch 5.3), has received special attention for two reasons. First, this implies that $\alpha_i$ can be interpreted as the (average) infectivity of $i$-individuals and $\beta_j$ as the (average) susceptibility of $j$-individuals. Secondly, the basic reproduction then has the simpler form $R_0=\sum_i\alpha_i\beta_i\pi_i$.

As for the homogeneous case only minor outbreaks are possible if $R_0\le 1$ whereas a large outbreak may occur if $R_0>1$. There is a related threshold limit theorem stating what the probability of a major outbreak is, and a central limit theorem for the final number of infected of the different types in case the epidemic takes off. The expressions are more involved as are the proofs giving the desired results, why we refer to Ball and Clancy (1993) for details.

\subsection{Mixing heterogeneity: household and network models}\label{mixing}
Perhaps even more important than allowing individuals to be different in terms of susceptibility and infectivity is to allow for non-uniform mixing, meaning that an individual has different (average) contacts rates with different individuals.
The multitype epidemic model of the previous section may include some such non-uniform mixing in the sense that contacts with different types of individuals may be different. However, an assumption of the contact rates between specific pairs of individuals in the multitype model is that they are all of order $1/n$. In many real-life situations individuals tend to have a few other individuals with which they mix at a much higher rate. Epidemic models taking such type of mixing behaviour into account are often referred to as two-level or multilevel mixing epidemic models. The two most common examples are \emph{household} epidemics and \emph{network} epidemics.

In a household epidemic individuals are grouped into small groups (households) and it is assumed that the contact rate between pairs of individuals of the same household equals $\lambda_H$ and contact rates between pairs of individuals of different households equals $\lambda_G/n$ (and all individuals have i.i.d.\ infectious periods $I$ with mean $1/\gamma$). Large population properties for this model was first rigorously analysed by Ball et al.\ (1997). There they show that it is necessary take into account the random distribution of the total number of infected in a household outbreak. During the initial phase of an outbreak, external infections are nearly always with not yet infected households, so by treating households as ''super-individuals'' the initial phase of the epidemic may be approximated by a multitype branching process, where the different types refer to household size and also the size of the household outbreak. The basic reproduction number (now having a more complicated definition but still having the value of 1 as its threshold) equals $R_0=(\lambda/\gamma)\mu_H$, where $\mu_H$ is the mean size of a household outbreak where the size of the selected household has a size biased distribution due to the fact that larger households are more likely to be externally infected as more individuals reside in them.

As before, major outbreaks are only possible if $R_0>1$, and when this happens Ball et al. (1997) show a law of large numbers for the number of households of different sizes and outbreak sizes using a balance equation. With some additional effort they also derive a central limit theorem for the final size in case of a major outbreak.

Another type of epidemic model having a local structure with much
higher (or all) contact rates between ''neighbours'' are called
network epidemic models (e.g.\ Newman, 2003). The network, an
undirected graph, specifies the social structure in the community upon
which an epidemic spreads. Quite often the network is assumed to be
random but having some pre-specified properties. These could for
example be the degree-distribution (the distribution of the number of
neighbours), the clustering coefficient (specifying how frequent
close cycles are present), and the degree correlation of
neighbours. Given such a network and a stochastic epidemic model
''on'' the network, it is of interest to see if a major outbreak can
occur, and if so what is its probability and what is the size of such an
outbreak. There are still many open problems to be solved in this area
(in particular when considering a dynamic network) but see for
example Britton et al.\ (2007) for some results. In general, the
social structure in the community are more influential for diseases
which are not highly infectious. One such class of diseases are
sexually transmitted infections -- in this case an edge between
two individuals in the network correspond a sexual relation.

\subsection{Models for endemicity}\label{endemicity}
The focus of the present paper has been SIR epidemic models for a
closed population, i.e.\ not allowing deaths or that new individuals
enter the community during the outbreak. This is of course an
approximation of real life but, if focus is on short time behaviour it
is sensible to make such an approximation. Some infectious diseases
are endemic in many countries, and a question of interest (e.g.\
Anderson and May, 1991) for such disease are to understand their
dynamics; to understand why certain diseases are endemic and others
not, and why a given disease may be endemic in one country but not in
another, and it is of course also of interest to determine which
preventive measure that has the potential of eradicating the
disease. We now give a brief outline to this problem area and refer
the reader the literature for more details; e.g.\ Anderson and May
(1991, pp 128) and N{\aa}sell (1999), and, when studying measles in
particular, e.g.\ Conlan and Grenfell (2007).
We now define the Markovian SIR epidemic model with demography (N{\aa}sell, 1999).

The population dynamics are very simple: new (susceptible)
individuals enter the population according to the time points of a
homogeneous Poisson process at rate $n\mu$, and each individuals lives
for an exponentially distributed time with mean $1/\mu$. In words,
there is a steady and constant influx of susceptibles at rate $\mu n$
and each individuals dies at rate $\mu$, so the population size will
fluctuate around the equilibrium state $n$, which hence can be
interpreted as ''population size''.

\begin{figure}[!h] \begin{center} \bf
\includegraphics[height=6cm, width=8cm]{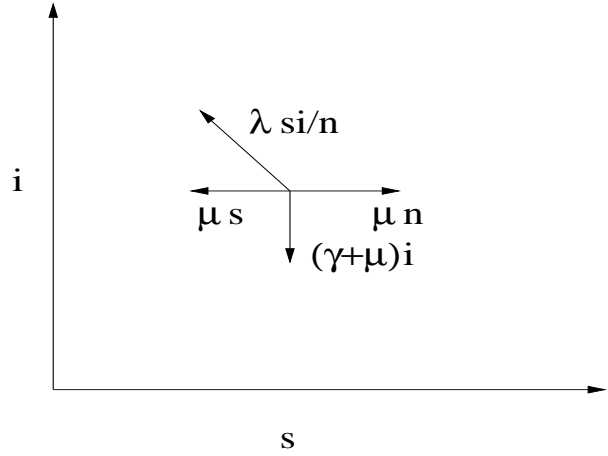}
\caption{\rm Jump rates for the SIR epidemic model with
  demography. All jumps are one unit of length except the jump to the
  north west in which $s$ decreases and $i$ increases by one unit (due to someone
  getting infected).}
 \label{fig:sir-dem}\end{center}
\end{figure}

The disease dynamics are just like for the Markovian version of the
standard stochastic epidemic model: an infectious individual has close
contacts with each other individual at rate $\lambda /n$: if such a
contact is with a susceptible individual this individual gets infected
and infectious immediately, and each
infectious individual recovers at rate $\gamma$ and becomes immune for
the rest of its life. Each individual dies
at rate $\mu$ irrespective of the infection status. The epidemic
process $\{ S(t), I(t), R(t); t\ge 0\}$ is initiated by $(S(0), I(0), R(0))=(s_0, i_0,r_0)$, where, as
before, we assume that $i_0>0$. In Figure \ref{fig:sir-dem} the various jump rates of
the model are given in the $(s,i)$-plane.
As before it is enough to keep track the number of
susceptibles and infectives, since recovered and immune individuals
play no further roll in the disease dynamics. A question of interest
is to study properties of this model when $n$
is large; we do this by first studying the corresponding deterministic
system defined by the differential equations:
\begin{align}
s'(t)&= \mu -\lambda s(t)i(t) - \mu s(t),\nonumber\\
i'(t)&= \lambda s(t)i(t)-\gamma i(t) - \mu i(t),\label{def-det-end}\\
r'(t)&= \gamma i(t) - \mu r(t).\nonumber
\end{align}
By equating all derivatives in (\ref{def-det-end}) to 0 we get the equilibrium state(s). It
is straightforward to show that $(s(t), i(t), r(t))=(1,0,0)$, the
disease free state, is always a point of equilibrium. The basic
reproduction number $R_0$, defined as the expected number of
individuals a typical individual infects during the early stages of an
outbreak, equals $\lambda /(\gamma + \mu)$ since now there are two
possible reasons for leaving the infectious state: recovering and dying.
Additional to the disease free equilibrium (which is stable if $R_0\le
1$ and unstable otherwise) there is another stable equilibrium whenever $R_0>1$:
\begin{equation}
(s(t), i(t), r(t))=(\hat s, \hat i, \hat r)=(\frac{1}{R_0},\ 
\frac{R_0-1}{R_0/\delta},\ 1-\frac{1/\delta + R_0-1}{R_0/\delta} ),\label{end-level}
\end{equation}
where $\delta=\mu/(\gamma + \mu)$ is the (small) average fraction of a
life that an individual is infectious. This state is called the
\emph{endemic level}.

For the stochastic counterpart it is possible to reach the state $I(t)=0$
from any other given state in finite time with positive probability. This observation
together with the fact that the state $I(t)=0$ is absorbing -- once
the disease disappears it can never return -- makes all
other states transient. As a consequence there is only one stationary
equilibrium: the disease free state in which all individuals
(fluctuating around the number $n$) are susceptible. However, when
$R_0>1$ and $n$ is large there is a drift of $(S(t), I(t), R(t))$
towards the endemic level $(n\hat s, n\hat i, n\hat r)$. So,
even though the disease eventually will go extinct, this may take a long time,
and in the mean time the epidemic process will fluctuate around the
endemic level. This type of behaviour is known as
\emph{quasi-stationarity} (e.g.\ van Doorn, 1991). Questions of interest are for
example, given the population size $n$ and the parameters $\lambda$,
$\gamma$ and $\mu$ (or equivalently $\mu$, $R_0$ and $\delta$) how
long will it take for the disease to go extinct? Or, phrased in
another way, how big must the population be for a given disease to persist
in the community, and which features of the disease are most
influential in determining this so called \emph{critical community
  size}? Another important question is to study effects of
introducing vaccination in the model.
We refer to for example N{\aa}sell (1999) and Andersson and Britton
(2000b) for stochastic methodological work in the area of endemic diseases.

\subsection{Other extensions}\label{other-ext}

We have focused on presenting results for a fairly
simple stochastic epidemic model; the reason being that even in simple
models results are far from trivial. Reality is of course much more
complicated. There are many features that affect the spread of
infections that have not been discussed in the paper - below we list
some such features, but the are also many others.

It
is well-known that seasonal effects play an important roll in disease
dynamics. One reason is that certain viruses or bacteria reproduce at
higher speed under certain seasons, e.g.\ influenza virus prefers cold
weather, but perhaps even more important is the change in social
behaviour over the year. The classical example is school semester and
the school start in September being the event that triggered measles
outbreaks in for example England prior to vaccination (e.g.\ Anderson and May, 1991, Ch 6). One of the most recent
illustrations of effects of schools is by Cauchemez et al.\ (2008),
where the influence of transmission within school is estimated by
comparing over-all transmission during the school semester with over-all transmission during holidays.

The model studied in the present paper lacks a spatial
component. Even though travel has increased dramatically the last
century there is of course still a spatial component in the spreading
of infections, both between and within countries. In larger
modelling/simulation studies of specific diseases, like the new
influenza pandemic, space is always taken into account (e.g.\ Yang et
al., 2009,
and Fraser et al., 2009). This is done by assuming
that the probability to infect a given individual decays with the
distance between the individuals. The main effect of introducing space into
the model is that the epidemic growth is slower: more linear than
exponential.

In the model it was assumed that the infection rate $\lambda$ was
constant during the infectious period. This rate can be thought of as
the product of two quantities: the rate of having contact with other
individuals and the probability of transmission upon such a
contact. In reality both of these quantities typically vary over the
infectious period. First there is often a latency period when an
infected individual hardly is infectious, then there is often a period
where the individual is highly infectious but still have few or no
symptoms. Eventually the individual has symptoms thus reducing the
social activity, and finally the infectiousness starts dropping down
to 0. This behaviour will clearly affect the dynamics of the disease
(a disease having a long latency period of course slows down the
epidemic growth), but it has been shown (Ball, 1983 and 1986) that final size
results of the present model still applies. The total infectivity of
an infectious period in the present model equals the random quantity
$\lambda I$. For a model in which the infection rate $t$ time units after
infection equals $C(t)\iota (t)$, where $C(t)$ and $\iota (t)$ are
stochastic processes describing social activity and infectiousness,
then the new model may be analysed with the present model if the
distribution of $\lambda I$ is replaced by the distribution of
$\int_0^\infty C(t)\iota (t)dt$. In particular, a latency period of
arbitrary length has no effect on the distribution of the final size.

In the section about modelling vaccination it was assumed that the
vaccine was perfect in the sense that it gave complete immunity to the
disease. In reality this is rarely the case. There are several more
realistic models for vaccine response in which the vaccine reduces
susceptibility, symptoms and/or infectivity in case of being infected,
where all these reductions are random (e.g.\ Halloran et al., 2003,
and Becker et al., 2006). For example, a
vaccine which reduces susceptibility by a factor $e$,
so the  average relative susceptibility is $(1-e)$ as compared with unvaccinated,
but has no effect on infectivity in case an infection occurs, has a new
(higher) critical vaccination coverage
$$
v_c=\frac{1}{e}\left(1-\frac{1}{R_0}\right) .
$$
This true if all individuals have the same reduction $e$ (so-called
\emph{leaky} vaccines) as well as if a fraction $e$ of the vaccinated become
completely immune and the rest are unaffected by the vaccine
(\emph{all-or-nothing} vaccines), or something in-between.
Having
models for vaccine efficacy is of course not enough for making conclusions in specific
situations -- it is equally important to estimate the various
\emph{vaccine efficacies}. This is most often done in clinical trials,
and it is usually harder to estimate reduction in infectivity because
this has to be done indirectly since actual infections are rarely
observed,  see for example Becker et al.\ (2006) and the forth coming book
by Halloran et al.\ (2010).

 \section*{Discussion}

Even when trying to include as many realistic features in a model as
possible there is a limit to how close a model can get to reality, and models can never completely predict what will happen in
a given situation. It is for example nearly impossible to predict how
people will adapt and change behaviour as a disease starts
spreading. Having said this, models can still be (and are!) very
useful as guidance for health professionals when deciding about
preventive measures aiming at reducing the spread of a disease.

Stochastic epidemic models, or minor modifications of them, can be
used also in other areas. A classic example is models for the spread
of rumours or knowledge (e.g.\ Daley and Gani, 1999, Ch 5). More
recently, the world wide web have several aspects resembling epidemic
models: for example computer viruses (which even use terminology from
infectious diseases) and spread of information.

The present paper only gives a short introduction to this rather big
research field. There are several monographs about mathematical models
for infectious disease spread: Anderson and May (1991) is probably the
most well-known, Diekmann and Heesterbeek (2000) has a higher
mathematical level, and  Keeling and Rohani (2008) also explicitly
considers disease spread among animals.
When it comes to stochastic epidemic models  there are for example the
classic book by Bailey (1975), and Andersson and Britton (2000a) who
also cover inference; a topic which Becker (1989) focuses on.

 \section*{Acknowledgements}
 I am grateful to the Swedish Research Council for financial support.


\begin{thebibliography}{99}



\bibitem{am91} Anderson R. M. and May R. M. (1991).
{\it Infectious diseases of humans; dynamic and control}.
Oxford: Oxford University Press.


\bibitem{ab00a} Andersson H. and Britton T. (2000a).
{\it Stochastic epidemic models and their statistical analysis}.
Springer Lecture Notes in Statistics. New York: Springer Verlag.

\bibitem{ab00b} Andersson H. and Britton T. (2000b): Stochastic
  epidemics in dynamic populations: quasi-stationarity and
  extinction. \emph{J. Math. Biol.}, {\bf 41}, 559-580.


\noindent von Bahr, B.\ and  Martin-L\"of, A.\ (1980): Threshold limit
theorems for some epidemic processes, {\it Adv. Appl. Prob., \bf 12},
319-349. 


\bibitem{b75} Bailey,  N.T.J. (1975).
{\it The Mathematical Theory of Infectious Diseases and its Applications}.
      London: Griffin.



\bibitem{b83} Ball, F.G.\ (1983). The threshold behaviour of epidemic
  models. {\it J. Appl. Prob. \bf 20}, 227-241.


\bibitem{ball86} Ball F. G. (1986). A unified approach to the distribution
of total size and total area under the trajectory of the infectives in
epidemic models. {\it Adv. Appl. Prob. \bf18}, 289-310.


\bibitem{bc93} Ball F. G., and Clancy D. (1993). The final size and severity
of a generalised stochastic multitype epidemic model. {\it Adv. Appl.
Prob. \bf 25}, 721-736.

\bibitem{bms97} Ball F. G., Mollison D. and Scalia-Tomba G. (1997).
Epidemics with two levels of mixing. {\it Ann. Appl. Prob. \bf 7}, 46-89.


\bibitem{B74} Barbour A. D. (1975): The duration of the closed stochastic
epidemic. {\it Biometrika, \bf 62}, 477-482.


\bibitem{ba49} Bartlett  M. S. (1949). Some evolutionary stochastic
processes. {\it J. Roy. Statist. Soc. B, \bf 11}, 211-229.


\bibitem{bec89} Becker N. G. (1989). {\it Analysis of Infectious Disease
 Data}. London: Chapman and Hall.


\bibitem{BBO06} Becker N.G., Britton T. and O'Neill P.D. (2006):
  Estimating vaccine effects from studies of outbreaks in household
  pairs. \emph{Statistics in Medicine}, {\bf 25}, 1079-1093.


\bibitem{db1760}
Bernoulli, D. (1760) Essai d'une nouvelle analyse de la mortalit\'e
caus\'ee par la petite v\'erole et des avantages de l'inoculation
pour la pr\'evenir, M\'em. Math. Phys. Acad. Roy. Sci., Paris (1760)
1--45.


\bibitem{bollobas}
Bollob\'as, B.\ (2001) \emph{Random Graphs}. 2nd ed., Cambridge Univ. Press,
Cambridge.



\bibitem{BJM07} Britton, T., Janson, S., Martin-Löf A. (2007): Graphs
  with specified degree distributions, simple epidemics and local
  vacination strategies. \emph{Adv. Appl. Prob.}, {\bf 39}, 922-948. 



\bibitem{C08}
Cauchemez, S., Valleron, A., Bo{\"e}lle, P., Flahault, A., and Ferguson, N.
  (2008). Estimating the impact of school closure on influenza transmission
  from sentinel data. \emph{Nature}, {\bf 452}, 750--754.


\bibitem{CG07} Conlan A.J.K., Grenfell B.T.\ (2007). Seasonality and the persistence
and invasion of measles. \emph{Proc. Roy. Soc Lond B} {\bf 274(1614)}: 1133-1141.


\bibitem{AC98} Cox, C. (1998). The delta method. In \emph{Encyclopedia
    of Biostatistics}, P Armitage and T. Colton (Eds). Wiley: Chichester.


\bibitem{DG99} Daley D. J. and Gani J. (1999). {\it Epidemic Modelling: an
introduction}. Cambridge University Press, Cambdridge.


\bibitem{DH00} Diekmann, O.\ and Heesterbeek, J.A.P.\ (2000): {\it
Mathematical Epidemiology
of Infectious Diseases: Model Building, Analysis and Interpretation.} Wiley, Chichester.


\bibitem{vD91} van Doorn, E.A. (1991). Quasi-stationary distributions
and convergence to quasi-stationarity of birth-death
processes. \emph{Adv. Appl. Prob.}, {\bf 23}, 683-700.


\bibitem{F} Fraser, C., Donnelly, C.A., Cauchemez, S., Hanage, W.P.,
  Van Kerkhove, M.D., Hollingsworth, T.D., Griffin, J., Baggaley,
  R.F., Jenkins, H.E., Lyons, E.J., et al. (2009). Pandemic Potential
  of a Strain of Influenza A (H1N1): Early Findings. \emph{Science},
  {\bf 324}, 1557-1561.


\bibitem{HJV05} Haccou, P., Jagers, P.\ and Vatutin, V. A.\ (2005).
\emph{ Branching Processes: Variation, Growth, and Extinction of
Populations}. Cambridge University Press, Cambridge.


\bibitem{HLS10} Halloran, M.E., Longini, Jr., I.M., Struchiner,
  C.J. (2010). \emph{Design and Analysis of Vaccine Studies}, Springer.


\bibitem{HPC03} Halloran, M.E., Pr\'eziosi, M-P., Chu, H. (2003)
Estimating vaccine efficacy from secondary attack rates. {\em
Journal of the American Statistical Assocociation}, {\bf
98}, 38-46.


\bibitem{J75} Jagers, P. (1975) {\it Branching Processes with
Biological Applications}. Wiley, London.


\bibitem{KR07}
Keeling, M. and Rohani, P. (2008).
{\em Modeling infectious diseases in humans and animals}.
Princeton University Press, Princeton.


\bibitem{K56} Kendall, D.G.\ (1956) Deterministic and stochastic
  epidemics in closed populations. \emph{Proc.\ Thirs Berkeley
    Symp.\ Math.\ Statist.\ \& Prob.}, {\bf 4}, 149-65.

\bibitem{KM27} Kermack, W. O. and McKendrick, A. G. (1927). A
contribution
to the mathematical theory of epidemics. {\it Proc. Roy.
Soc. Lond. A, \bf 115}, 700-721.

\bibitem{m77} Mollison D. (1977). Spatial contact models for ecological and epidemic spread (with Discussion). {\it J. Roy. Statist. Soc. B
\bf39}, 283-326.

\bibitem{n03} Newman, M.E.J.\ (2003): The structure and function of complex
networks, \emph{SIAM Rev.} \textbf{45}, 167-256.


\bibitem{N99} N{\aa}sell, I. (1999): On the time to extinction
in recurrent epidemics. {\it J. Roy. Statist.
Soc. B.\bf 61}, 309-330.


\bibitem{PL90} Picard, P. and Lef\`{e}vre, C.\ (1990). A unified analysis of the final size
and severity distribution in collective Reed-Frost epidemic processes.
\emph{ Adv. Appl. Prob.}, {\bf 22}, 269--294.


\bibitem{R11} Ross, R (1911). \emph{The prevention of malaria, 2nd
    ed.}. London: Murray.


\bibitem{S85} Scalia-Tomba, G.\ (1985). Asymptotic final size
  distribution for some chain-binomial processes. \emph{Adv.\ Appl.\
    Prob.}, {\bf 17}, 477-495.


\bibitem{S90} Scalia-Tomba, G.\ (1990). On the asymptotic final size
  distribution of epidemics in heterogeneous populations.  In in
  Stochastic Processes in Epidemic Theory, Gabriel, J.-P., 
Lef\`evre, C.,  Picard, P. (Eds.). {\it Springer Lecture
Notes in Biomath.\ \bf 86}: 189-196.


\bibitem{S83} Sellke, T.\ (1983). On the asymptotic distribution of the size of a stochastic
epidemic. \emph{J.~Appl.~Prob.} {\bf 20}, 390-394.


\bibitem{Y09} Yang, Y., Sugimoto, J.D., Halloran, M.E., Basta, N.E., Chao, D.L.,
  Matrajt, L., Potter, G., Kenah, E., Longini, I.M. Jr. (2009) The
  Transmissibility and Control of Pandemic Influenza A (H1N1) Virus. \emph{Science}. [DOI: 10.1126/science.1177373] .


\end{thebibliography}
\end{document}